\documentclass[10pt,twoside,reqno]{amsart}
\usepackage{amsmath,amssymb,amsfonts,amsthm,epsfig,xspace}
\usepackage{ecphead}
\newcommand{\old}[1]{}
\old{
\setlength{\vsize}{9in}
\setlength{\hsize}{6.5in}
\setlength{\textwidth}{6.5in}
\setlength{\textheight}{9in}
\setlength{\oddsidemargin}{0in}
\setlength{\evensidemargin}{0in}
\setlength{\topmargin}{0in}
}

\newcommand{\eps}{\varepsilon}
\newcommand{\Z}{{\mathbb Z}}
\newcommand{\N}{{\mathbb N}}
\newcommand{\E}{{\mathbb E}}

\newtheorem{theorem}{Theorem}
\newtheorem{lemma}[theorem]{Lemma}
\newtheorem{proposition}[theorem]{Proposition}

\newcommand{\lref}[1]{Lemma~\ref{lem:#1}}

\newcommand{\eref}[1]{Equation~\eqref{eqn:#1}}

\newcommand{\SRW}{\operatorname{SRW}}

\begin{document}
\input{bw-head.sty}
\old{
\title[Excited random walk]{\vspace*{-50pt}
 Excited Random Walk}
\author[Benjamini \& Wilson]{Itai Benjamini \and David B. Wilson}
\address{Weizmann Institute\\Rehovot 76100, Israel}
\address{Microsoft Research\\One Microsoft Way\\Redmond, WA 98052, U.S.A.}


\vspace*{-12pt}
\begin{abstract}
A random walk on $\Z^d$ is excited if the first time it visits a
vertex there is a bias in one direction, but on subsequent visits
to that vertex the walker picks a neighbor uniformly at random.
We show that excited random walk on $\Z^d$ is transient iff $d>1$.
\end{abstract}

\maketitle
}

\section{Excited Random Walk}

A random walk on $\Z^d$ is excited (with bias $\eps/d$) if the first time
it visits a
vertex it steps right with probability $(1+\eps)/(2d)$ ($\eps>0$), left with
probability $(1-\eps)/(2d)$, and in other directions with probability
$1/(2d)$, while on subsequent visits to that vertex the walker picks a
neighbor uniformly at random.  This model was studied heavily in the
framework of perturbing $1$-dimensional Brownian motion, see for
instance \cite{D2,PW} and reference therein.  Excited random walk
falls into the notorious wide category of self-interacting random
walks, such as reinforced random walk, or self-avoiding walks. These
models are difficult to analyze in general. The reader should consult
\cite{D,L,TW,S,ABV}, and especially the survey paper \cite{pemantle01}
for examples. Simple coupling and an additional neat observation allow
us to prove that excited random walk is recurrent only in dimension $1$.
The proof uses and studies a special set of points (``tan points'') for the
simple random walk.

\section{Recurrence in $\Z^1$}

It is already known that excited random walk in $\Z^1$ is
recurrent, indeed, a great deal more is known about it \cite{D3}.
But for the reader's convenience we provide a short proof.

On the first visit to a vertex there is probability $p>1/2$ of going
right and $1-p$ of going left, while on subsequent visits the
probabilities are $1/2$.  Suppose that the walker is at $x>0$ for the
first time, and that all vertices between $0$ and $x$ have been
visited.  The probability that the walker goes to $x+1$ before going
to $0$ is $p+(1-p)(1-2/(x+1)) = 1 - 2(1-p)/(x+1)$.  Multiplying over
the $x$'s, we see that the random walk returns to $0$ with probability $1$.

\section{Transience in $\Z^2$}

The simple random walk (SRW) in $\Z^2$ visits about $n/\log n$
points by time $n$, and if the excited random walk (ERW) gets
pushed to the right $n/\log n$ times, it would very quickly depart
its start location and never return.  But it is not clear what
effect that the perturbations have on the number of visited
points, and it is not obvious that the excited random walk will
visit $n/\log n$ distinct points by time $n$.

To lower bound the number of points that the excited random walk
visits, we couple it with the SRW in the straightforward way, and
count the number of ``tan points'' visited by the SRW.  We define
the coupling as follows: if the SRW moves up, down, or right, then so does the
ERW. If the SRW moves left, then the ERW moves left if it is at a
previously visited point, and if the ERW is at a new point, it moves
either left or right with suitable probabilities.  At all times,
the $y$-coordinates of the SRW and ERW are identical.

To explain the concept of a ``tan point'', we imagine that the
simple random walker leaves behind an opaque trail, and that the
sun is shining from infinitely far away in the positive
$x$-direction.  If the SRW visits a point $(x,y)$ such that no
point $(x',y)$ with $x'>x$ has been visited, then the sun shines
upon $(x,y)$, and this point becomes tanned.  Formally, we define a
\textit{tan point\/} for the SRW to be a vertex $(x,y)$ that is visited
by the SRW before any point of the form $(x',y)$ with $x'>x$.
If the sun shines upon the simple random walker the first time it
is at $(x,y)$, it is straightforward to check that ERW is at a new
point.  We will show that with high probability there are many tan
points (so the ERW visits many new points), and that this implies
that the ERW is transient.

The probability that a point $(x,y)$ will be tan follows from some
enumerative work of Bousquet-M\'elou and Schaeffer on random walks in
the slit plane \cite{BMS}.
\begin{lemma}
Let $r$ and $\theta$ be the polar coordinates of the point $(x,y)$,
i.e.\ $r\geq 0$, $0\leq\theta<2\pi$, $x=r\cos\theta$, and $y=r\sin\theta$.
Then
\begin{equation}\label{eqn:pr-tan}
\Pr[\text{$(x,y)$ is tan}] = (1+o(1))\sqrt{\frac{1+\sqrt{2}}{2\pi}} \frac{\sin(\theta/2)}{\sqrt{r}},
\end{equation}
where the $o(1)$ term goes to $0$ as $r$ tends to $\infty$.
\end{lemma}
This equation does not explicitly appear in \cite{BMS}, but all the
real work that goes into proving it is in \cite{BMS}.  In the interest
of completeness, we explain how this equation follows from explicit
results in \cite{BMS}:
\begin{proof}
Let $a_n$ be the number of walks of length $n$ that start from $(0,0)$,
and avoid the nonnegative real axis at all subsequent times, and let
$p_{x,y,n}$ denote the probability that a random such walk ends at the
point $(x,y)$.  By reversibility of the random walks,
$$\Pr\left[\begin{aligned}\text{SRW started from the point $(x,y)$ first hits the}\\ \text{nonnegative real axis at the point $(0,0)$ and at time $n$}\end{aligned}\right] = \frac{a_n}{4^n} \times p_{x,y,n}.$$
Thus $\Pr[\text{$(x,y)$ is tan}]=\sum_{n=0}^\infty a_n/4^n \times
p_{x,y,n}$.  Theorem~1 of \cite{BMS} gives
$$\frac{a_n}{4^n}=(1+o(1))\frac{\sqrt{1+\sqrt{2}}}{2\Gamma(3/4)}
n^{-1/4}.$$ Theorem~21 of \cite{BMS} considers the endpoint
$(X_n,Y_n)$ of a random walk started from $(0,0)$ which avoids the
nonnegative real axis, and gives the limiting distribution of
$(X_n/\sqrt{n},Y_n/\sqrt{n})$.  This limiting distribution morally
determines the asymptotics of $p_{x,y,n}$ --- a \textit{local\/} limit
theorem \textit{would\/} determine the asymptotics --- and the authors
prove a local limit theorem for $Y_n/\sqrt{n}$ but not $X_n/\sqrt{n}$
let alone the joint distribution $(X_n/\sqrt{n},Y_n/\sqrt{n})$.
However, since the ordinary random walk has a local limit theorem (on
vertices such that $x+y\equiv n \mod 2$), one can take the limiting
distribution of $X_{(1-\eps)n},Y_{(1-\eps)n}$ and then run the walk
another $\eps n$ steps; upon sending $\eps$ to $0$ sufficiently
slowly, one can obtain a local limit theorem version of Theorem~21 of
\cite{BMS}:
$$p_{x,y,n} = (1+o(1))\frac{2}{\Gamma(1/4)} \frac{r^{1/2}}{n^{5/4}} e^{-r^2/n} \sin(\theta/2) \times \begin{cases}2 & x+y\equiv n \mod 2 \\ 0 & x+y\not\equiv n \mod 2\end{cases}$$
when $r=\Theta(\sqrt{n})$ and $\theta$ is bounded away from $0$ and $2\pi$.
Thus we obtain
\begin{align*}
 \sum_{n=\Theta(r^2)} \frac{a_n}{4^n} \times p_{x,y,n}
 &= (1+o(1)) \int_0^\infty \frac{\sqrt{1+\sqrt{2}}}{2\Gamma(3/4)} n^{-1/4} \times \frac{2}{\Gamma(1/4)} \frac{r^{1/2}}{n^{5/4}} e^{-r^2/n} \sin(\theta/2) \, dn\\
 &= (1+o(1)) \sqrt{\frac{1+\sqrt{2}}{2\pi}} \frac{\sin(\theta/2)}{\sqrt{r}}.
\end{align*}

Next we check that the terms when $n\ll r^2$ or $n\gg r^2$ contribute negligibly.  We may bound $\frac{a_n}{4^n} \times p_{x,y,n}$ by the probability that the walk survives the first $n/2$ steps ($O(1/n^{1/4})$) times the probability that ordinary SRW for the remaining $n/2$ steps ends at the point $(x,y)$ ($O(1/n)$).  Thus 
\begin{align*}
\sum_{n\gg r^2} \frac{a_n}{4^n} \times p_{x,y,n} &\leq
\sum_{n\gg r^2} O(1/n^{5/4}) \ll 1/r^{1/2}
\end{align*}
When $n\ll r^2$ we bound $\frac{a_n}{4^n} \times p_{x,y,n}$ by the probability that the walk makes it out to radius $r/2$ without hitting the line (which is $O(1/r^{1/2})$ by \cite[Eqn~2.40]{L}) times the probability that the walk ends up at $(x,y)$ (distance $s\geq r/2$ away) at the end of the remaining $m<n$ steps.  This latter probability is at most $\frac{1+o(1)}{\pi m} e^{-s^2/(2m)}$, and assuming $n\leq r^2/8$, it is upper bounded by $\frac{1+o(1)}{\pi n} e^{-r^2/(8n)}$.  Thus
\begin{align*}
\sum_{n\ll r^2} \frac{a_n}{4^n} \times p_{x,y,n} &\leq
O(1/r^{1/2}) \times \sum_{n\ll r^2} O(n^{-1} e^{-r^2/(8n)}) \ll 1/r^{1/2}\\
\end{align*}
Hence when $\theta$ is bounded away from $0$ and $2\pi$ and $r$ is large, the terms when $n\ll r^2$ or $n\gg r^2$ contribute negligibly, so the formula follows in this case.  But the formula in the limiting case when $\theta$ approaches $0$ or $2\pi$ follows from the case when $\theta$ is bounded away from $0$ and $2\pi$, so the formula is valid simply when $r$ is large enough.
\end{proof}

We will use the notation introduced by Knuth where $\Theta(f)$ denotes
an expression which is upper-bounded by $C \times f$ and lower-bounded
by $c \times f$, where $c$ and $C$ are positive constants.  (By contrast,
$O(f)$ denotes an expression for which there is an upper bound of $C\times f$,
but not necessarily any lower bound.)
Using this notation, we may crudely approximate \eref{pr-tan} with
\begin{equation}\label{eqn:crude-tan}
\Pr[\text{$(x,y)$ is tan}] = \begin{cases}\Theta(1/r^{1/2}) & \text{when $x\geq 0$}\\ \Theta(|y|/r^{3/2}) & \text{when $x\leq 0$.} \end{cases}
\end{equation}

\eref{crude-tan} is useful, but we need two
modifications. For our purposes, it would be better to have the
probability that a point is tan and that it is reached by time
$n$.  Then we would have the expected number of tan points by time
$n$, which would lower-bound the expected number of times that the
ERW is pushed to the right.  But a lower bound on the expected
number is not quite what we need to prove transience of the ERW;
what we'd really like to know is that with very high probability,
the number of tan points by time $n$ is large.

In order to get the ``with very high probability'' part of the
statement, it would be convenient to be working with independent
events.  To get this independence, we divide the plane into bands of
height $h=h(n)$ to be determined later, and we will focus on every
other band, say the even ones.  When the SRW first arrives at an even band,
let us count the tan points within the band that are encountered before
the random walk reaches a different even band.  These counts are independent
for the different bands.  After $n$ steps it is
likely that order $\sqrt{n}/h$ even bands have been crossed, so it is
likely that the number of tan points dominates a sum of
$\Theta(\sqrt{n}/h)$ independent random variables.

\begin{lemma}
\label{lem:tan-band}
Consider a band of height $h$, i.e., $\Z\times[y_0,y_0+h-1]$.  After
the SRW first reaches this band, with probability $\Theta(1)$ the SRW
hits $\Theta(h^{3/2})$ tan points within the band before leaving the
enclosing band $\Z\times[y_0-h,y_0+2 h-1]$.
\end{lemma}

For the $\Theta(1)$ part of this lemma, the following proposition is useful.
\newcommand{\half}{{\textstyle\frac12}}
\newcommand{\quart}{{\textstyle\frac14}}
\begin{proposition}
\label{prp:chance-big}
If $X$ is a real-valued random variable and $\E[X]\geq 0$, then
 $$\Pr\left[X\geq \half \E[X]\right] \geq \frac{\E[X]^2}{4\E[X^2]} .$$
\end{proposition}
This proposition is essentially exercise 1.3.8 of Durrett \cite{Durrett}.
\old{
\begin{proof}
\begin{align*}
\E[X] &= \E\left[X|X\geq \half \E[X]\right] \Pr\left[X\geq \half \E[X]\right] + \E\left[X|X<\half \E[X]\right] \Pr\left[X<\half \E[X]\right] \\
\E[X] &\leq \E\left[X|X\geq \half \E[X]\right] \Pr\left[X\geq \half \E[X]\right] + \half\E[X]\\
\quart\E[X]^2 &\leq \E\left[X|X\geq \half \E[X]\right]^2 \Pr\left[X\geq \half \E[X]\right]^2 \\
 &\leq \E\left[X^2|X\geq \half \E[X]\right] \Pr\left[X\geq \half \E[X]\right]^2 \\
 &\leq \E[X^2] \Pr\left[X\geq \half \E[X]\right] \qedhere
\end{align*}
\end{proof}
}

\begin{proof}[Proof of \lref{tan-band}]
Consider a point $(x,y)$ within the band, and the ray $[x,\infty)\times y$ with $(x,y)$ at its tip.
Consider the largest circle contained within the enclosing band and centered at the tip, and also a small enough disk centered at the tip.  From \eref{crude-tan} it follows that we can take the ratios of the two radii to be $\Theta(1)$ and have the property that for any point $p$ in the left half of the small disk and any point $q$ on the outer circle,
  $$\Pr[\text{$\SRW_p$ hits ray at tip}] \geq 2 \Pr[\text{$\SRW_q$ hits ray at tip}],$$
where $\SRW_p$ denotes the simple random walk started at point $p$, and by ``hits ray at tip'' we mean that the first time the walker hits the ray $(x+,y)$ is at the tip $(x,y)$.  Now
\newcommand{\event}{\Pr\left[\begin{aligned}\text{$\SRW_p$ hits ray at tip}\\ \text{before leaving circle}\end{aligned}\right]}
\begin{align*}
\Pr[\text{$\SRW_p$ hits ray at tip}] &= \event + \Pr\left[\begin{aligned}\text{$\SRW_p$ leaves circle and}\\ \text{then hits ray at tip}\end{aligned}\right]\\
\Pr[\text{$\SRW_p$ hits ray at tip}] &\leq \event + \max_q \Pr[\text{$\SRW_q$ hits ray at tip}]\\
\Pr[\text{$\SRW_p$ hits ray at tip}] &\leq \event + \frac12\Pr[\text{$\SRW_p$ hits ray at tip}]\\
\event &\geq \frac12\Pr[\text{$\SRW_p$ hits ray at tip}] = \Theta(1/\operatorname{dist}(p,\text{tip})^{1/2})
\end{align*}
The point where the random walk first enters the band will lie within the left half of the small disk surrounding $\Theta(h^2)$ such points $(x,y)$.  In fact,
there are $\Theta(h^2)$ such points $(x,y)$ within radius $h$ of where the SRW first hits the band.  Thus
$$\E\left[\begin{aligned}\text{\# tan points within band before SRW departs enclosing band}\\ \text{and within radius $h$ from where SRW arrives in band}\end{aligned}\right]
\geq \Theta(h^2)\times \Theta(1/h^{1/2}) = \Theta(h^{3/2}). $$
Next we need a second moment estimate:
\begin{align*}
\E&\left[\left(\begin{aligned}\text{\# tan points within band before SRW departs enclosing band}\\ \text{and within radius $h$ from where SRW arrives in band}\end{aligned}\right)^2\right] \\
  &\leq
\E\left[\left(\begin{aligned}\text{\# tan points within radius $h$ from where SRW arrives in band}\end{aligned}\right)^2\right] \\
  &= \sum_{\text{$x_1$,$y_1$,$x_2$,$y_2$ within radius $h$}}
   \Pr[\text{$(x_1,y_1)$ and $(x_2,y_2)$ tan}]\\
  &\leq 2 \sum_{\text{$x_1$,$y_1$,$x_2$,$y_2$ within radius $h$}}
   \Pr[\text{$(x_1,y_1)$ tan and $(x_2,y_2)$ tan after $(x_1,y_1)$}]\\
\intertext{where by ``after'' we include the possibility that $(x_2,y_2) = (x_1,y_1)$}
  &= 2 \sum \begin{aligned}[t] &\Pr[\text{$(x_1,y_1)$ tan}] \times
            \Pr[\text{$(x_2,y_2)$ after $(x_1,y_1)$} |  \text{$(x_1,y_1)$ tan}] \times\\
            &\Pr[\text{$(x_2,y_2)$ tan} | \text{$(x_2,y_2)$ after $(x_1,y_1)$ and $(x_1,y_1)$ tan}] \end{aligned}\\
  &\leq 2 \sum \Pr[\text{$(x_1,y_1)$ tan}] \times \Pr\left[\begin{aligned}&\text{after $(x_1,y_1)$, $(x_2,y_2)$ visited before $(x_2+,y_2)$} |\\& \text{$(x_2,y_2)$ after $(x_1,y_1)$ and $(x_1,y_1)$ tan}\end{aligned}\right]\\
  &= 2 \sum \Pr[\text{$(x_1,y_1)$ tan}] \Pr[\text{$(x_2-x_1,y_2-y_1)$ tan in $\SRW_0$}]\\
  &= \Theta(h^3)
\end{align*}

Combining these estimates with Proposition~\ref{prp:chance-big}, we
see that with at least $\Theta(1)$ probability there are at least
$\Theta(h^{3/2})$ tan points in the band before the SRW departs the
enclosing band.
\end{proof}

\begin{theorem}
With probability $1$, for all but finitely $n\in\N$, the excited
random walk has drifted right by a distance of at least
$\Theta(n^{3/4} / \log^{5/4} n)$ at time $n$.  In particular, it
is transient.
\end{theorem}

\begin{proof}
Say that the SRW deals with a band if it reaches that band and then
reaches a different band of the same parity.  Suppose that the random
walk starts from an odd band in the middle of a group of $4 k+1$ bands
of height $h$.  Then the probability that it fails to leave the group
of bands after $n=(k h)^2 t$ steps is exponentially small in $t$.  We
will optimize $k$ and $t$ later, but we will take $t$ to be large so
that with high probability the walk leaves the group of $4 k+1$ bands,
and in particular deals with at least $k$ even bands.

Rather than run the random walk for exactly $n$ steps, let us run it
until it deals with $k$ even bands.  Then the number of early tan
points in the different even bands are independent of one another, and
each one has a $\Theta(1)$ chance of being at least $\Theta(h^{3/2})$.
Except with probability $\exp(-\Theta(k))$, the number of tan points
in the $k$ even bands will be $\Theta(k h^{3/2})$.

Since there is only a $\exp(-\Theta(t))$ chance that the walk has not
dealt with $k$ even bands by time $n$, we find that, except with
probability $\exp(-\Theta(t))+\exp(-\Theta(k))$, there are
$\Theta(k h^{3/2})$ tan points by time $n$.  To optimize our
parameters we take $t=k$, and then we have $n = k^3 h^2$.

Next we consider the location of the perturbed random walk at time
$n$.  Typically the random walk diffuses by
$\Theta(\sqrt{n})=\Theta(k^{3/2} h)$, and drifts right by at least
$\Theta(k h^{3/2})$.  The probability that it diffuses by more than
$k^2 h$ is $\exp(-\Theta(k))$, and the probability that it drifts
less than $\Theta(k h^{3/2})$ is $<\exp(-\Theta(k))$.  We take
$k=\Theta(\log n)$, so that the probability of a bad event is
$<1/n^2$, which gives us $h=\Theta(n^{1/2}/\log^{3/2} n)$.  Except with
probability $<1/n^2$, the excited random walk has drifted right by
at least $\Theta(n^{3/4} / \log^{5/4}n)$.  In particular this
event fails only finitely often, so the random walk is transient.
\end{proof}

\section{Transience in $\Z^d$, $d>2$}

As with ERW in $\Z^2$, in $\Z^d$ we can couple the ERW with SRW, and
then (in continuous time) couple the SRW in $\Z^d$ with the SRW in
$\Z^2$.  For each tan point of the SRW in $\Z^2$ there is a tan point
of the SRW in $\Z^d$, so as before the ERW is transient.

\section{Speed}
We have seen that the excited random walk on $\Z^d$ is transient for
$d\geq 2$, but does it have positive speed?
\begin{theorem}
Let $X_n$ denote the $x$-coordinate at time $n$ of the excited random walk on $\Z^d$ with bias $\eps/d$.  If $d\geq 4$, then almost surely $\liminf_{n\rightarrow\infty} X_n/n \geq 0.659 \eps/d$; in particular the speed is positive.
\end{theorem}
\begin{proof}
Project down the $x$-coordinate of the ERW and $d-4$ additional coordinates, and consider the resulting SRW on $\Z^3$.  Let $R_n$ be the range of the SRW on $\Z^3$ by time $n$, i.e.\ the number of points visited by time $n$.  Since the SRW is transient, $\E[R_n]/n\rightarrow c$ where $c$ is the escape probability of the SRW.  Glasser and Zucker \cite{GZ} (see also \cite{DS}) determined this escape probability $c$ to be $$c=\frac{32\pi^3}{\sqrt{6}\Gamma(1/24)\Gamma(5/24)\Gamma(7/24)\Gamma(11/24)}=0.65946\ldots.$$  For our purposes it is not enough to know $\E[R_n]$, what we need is the strong law of large numbers for $R_n$ that was proved by Dvoretzky and Erd\H os \cite{DE}: a.s.\ $R_n/n\rightarrow c$ (see also \cite{BK} for even stronger results on $R_n$).  Thus for any $\delta>0$, a.s.\ there are only finitely many $n$ for which the ERW has not had $(c-\delta) n$ pushes to the right by time $n$.  The theorem then follows from the ordinary strong law of large numbers.
\end{proof}

It seems intuitive that excited random walk in $\Z^3$ also has
positive speed, but we do not see a proof.  Excited random walk in
$\Z^2$ is more delicate, and it is not clear even at an intuitive
level whether or not the speed is positive, though we believe that by
time $n$ it has traveled distance at least $\Theta(n/\log n)$.

\section*{Acknowledgements}
We are grateful to Oded Schramm for useful discussions, and we thank
Gabor Pete and the referee for their comments on an earlier version of
this article.  The research leading to this article was conducted
while the first author was visiting Microsoft.

\bibliographystyle{plainnat}

\end{document}